\documentclass[10pt]{amsart} % generic ams journal style

\newtheorem{theo}{Theorem}

\newtheorem{prop}{Proposition}

\newtheorem{coro}{Corollary}

\newtheorem{lemm}{Lemma}

\theoremstyle{definition}

\theoremstyle{remark}

\newtheorem{rema}{Remark}[section]

\newcommand{\LW}{\text{LambertW}}

\newcommand{\Exp}{\mbox{Exp}}

\begin{document}

\title[Integration]{An infinite-dimensional manifold structure for analytic Lie pseudogroups of infinite type}

\author{ Niky Kamran }

\author{ Thierry Robart }

\address{Department of Mathematics, Mc Gill University, Montreal, QC H3A 2K6, Canada}
\email{nkamran@math.mcgill.ca }
\address{Department of Mathematics, Howard University, Washington, D.C., 20059, U.S.A.}
\email{trobart@fac.howard.edu}

%\subjclass{Primary 47A15; Secondary 46A32, 47D20}

%\keywords{geometry, partial differential equations, contact}

\date{August 8, 2003.}

\dedicatory{}

%%% ----------------------------------------------------------------------

\begin{abstract}

We construct a infinite-dimensional manifold structure adapted to
analytic Lie pseudogroups of infinite type. More precisely, we
prove that any isotropy subgroup of an analytic Lie pseudogroup of
infinite type is a regular infinite-dimensional Lie group,
modelled on a locally convex strict inductive limit of Banach
spaces. This is an infinite-dimensional generalization to the case
of Lie pseudogroups of the classical second fundamental theorem of
Lie.

\end{abstract}

%%% ----------------------------------------------------------------------

\maketitle

%%% ----------------------------------------------------------------------

\section{Introduction}
Our objective in this paper is to prove an integration theorem for
the class of infinite-dimensional Lie algebras of analytic vector
fields corresponding to infinitesimal actions of Lie pseudogroups
of infinite type. More precisely, we shall show that any isotropy
subgroup of an analytic Lie pseudogroup of infinite type is a
regular infinite-dimensional Lie group in the sense of Milnor
\cite{Mil83}, modelled on a locally convex strict inductive limit
of Banach spaces. This is an infinite-dimensional generalization
to the case of Lie pseudogroups of infinite type of the classical
second fundamental theorem of Lie for finite-dimensional local Lie
algebras of vector fields. The main result of our paper represents
the conclusion of a research program that we have been involved in
over the past several years \cite{KaR97-1},\cite{KaR97-2},
\cite{RK97},\cite{KaR00},\cite{KaR01}, the aim of which has been
to construct a natural infinite-dimensional manifold structure for
the parameter space of analytic Lie pseudogroups of infinite type.

Recall that an analytic Lie pseudogroup $\Gamma^{\omega}$ of
transformations of an analytic manifold $M$ is a sub-pseudogroup
of the pseudogroup of analytic local diffeomorphisms of $M$, where
the elements of $\Gamma^{\omega}$ are the analytic solutions of a
system $\mathcal{S}$ of differential equations which is involutive
in the sense of Cartan-K\"ahler theory. Lie pseudogroups of
infinite type correspond to the case in which the solutions of
$\mathcal{S}$ are parametrized by arbitrary functions, while Lie
pseudogroups of finite type correspond to the case where the
solutions of $\mathcal{S}$ are parametrized by arbitrary
constants. A simple example of an analytic Lie pseudogroup of
infinite type is given by the set of conformal local
diffeomorphisms of the plane, with $\mathcal{S}$ being the
Cauchy-Riemann equations. The standard local action of
$PSL(2,\mathbb{R})$ on the real line given by fractional linear
transformations defines a Lie pseudogroup of finite type, with
$\mathcal{S}$ being the third-order ordinary differential equation
expressing the vanishing of the Schwarzian derivative. The Lie
pseudogroups of finite type give rise to finite-dimensional local
Lie groups \cite{Olv96}. The modern theory of finite-dimensional
Lie groups can be viewed as a theory of global parameter spaces
for Lie pseudogroups of finite type.

The classical theory of Lie pseudogroups was founded by Lie and
significantly developed by Elie Cartan, who formulated a general
structure theory encompassing both Lie pseudogroups of finite and
infinite type. Cartan's approach is based on the construction of
$\Gamma^{\omega}$-invariant differential forms which generalize
the Maurer-Cartan forms to the case of Lie pseudogroups of
infinite type, and satisfy generalized Maurer-Cartan structure
equations \cite{Ca04}, \cite{Ca37}. These invariant differential
forms can be computed explicitly from the differential system
$\mathcal{S}$ defining the Lie pseudogroup \cite{Ca37}, and lead
to many important geometric applications,
\cite{OlvPo03-1},\cite{OlvPo03-2}. The notion of an involutive
G-structure provides a global geometric framework for Cartan's
theory which is well adapted to the study of many classes of Lie
pseudogroups, \cite{Ch53}. The Lie pseudogroups arise in this
context as local analytic automorphisms of various reductions of
the linear frame bundle of the manifold $M$, and the Maurer-Cartan
forms are given by the vector-valued canonical one-form to the
G-structure. The Cartan-K\"ahler Theorem and the formulation of
the involutivity of $\mathcal{S}$ in terms of the characters of
the Lie-algebraic tableaux associated to the generalized
Maurer-Cartan equations are cornerstones of this theory.

The works of Spencer and of his school \cite{KS72} provide an
alternative approach to the study of Lie pseudogroups and their
deformations, in which the question of formal integrability is
expressed in purely cohomological terms. The fundamental papers by
Malgrange ~\cite{Mal72} on Lie equations give a clear and
beautiful account of the main integration theorem for formally
integrable differential systems which is central to the work of
Spencer, and of the links between this result and the
Cartan-K\"ahler Theorem. Malgrange's account of the
Cartan-K\"ahler Theorem will play an important role in our paper.

Every analytic Lie pseudogroup gives rise to an
infinite-dimensional filtered Lie algebra of vector fields ( more
precisely a sheaf of filtered Lie algebra germs ). These Lie
algebras and the geometries for which they act as infinitesimal
automorphisms have been studied in depth in a number of important
papers, particularly in the transitive case, \cite{guillemin 64},
\cite{SiSe65}. It is this local correspondence which is the
starting point of our work. In \cite{KaR01}, we proved an
integration theorem for local Lie algebras of analytic vector
fields corresponding to {\em transitive flat} Lie pseudogroups of
infinite type. More precisely, we proved that any isotropy
subgroup of such a Lie pseudogroup is a regular
infinite-dimensional Lie group in the sense of Milnor
\cite{Mil83}, modelled on a locally convex strict inductive limit
of Banach spaces. In the present paper, we dispense with the
hypothesis of flatness - a hypothesis which roughly amounts to
requiring that the defining equations of the Lie pseudogroup have
constant coefficients - and prove an integration theorem for the
filtered Lie algebras associated to isotropy subgroups of general
analytic Lie pseudogroups. The chart we construct in the non-flat
case has a structure similar to the the one we had obtained for
the flat case in \cite{KaR01}, in the sense that it is given by a
convergent infinite product of exponentials. However the specific
nature of the chart and of the proof of convergence are
significantly different in the non-flat case, both at the
conceptual and technical levels, and considerably more involved. A
key ingredient in the proof is the existence of a {\em bounded
filtered} basis for the underlying Lie algebra of vector fields.
The existence of this basis follows from the basic assumption of
involutivity of the underlying differential system and from an
important estimate established by Malgrange \cite{Mal72} in the
course of his proof of the Cartan-K\"ahler Theorem. We then derive
the estimates which are necessary to show the convergence of the
infinite product of exponentials adapted to this bounded filtered
basis. Our differentiable structure is thus similar to the one
obtained by Leslie in his fundamental work on the
infinite-dimensional manifold structure of the groups of analytic
diffeomorphisms of a compact analytic manifold, \cite{L82}.

\section{Notation and statement of the main theorem}

\subsection{Notation}

We consider $\mathbb{R}^{\nu}$ endowed with standard coordinates
$x=(x^1,\ldots,x^{\nu})$. For any $\nu$-multiplet
$\alpha=(\alpha_1,\ldots,\alpha_{\nu})$ of non-negative integers
we set \[\mid\alpha\mid=\alpha_1+\cdots+\alpha_{\nu},\quad
\mid\alpha\mid'=\mid\alpha\mid-1,\quad  \alpha !=\alpha_1
!\cdots\alpha_{\nu}!,\]
 and $$x^{\alpha}=(x^1)^{\alpha_1}\cdots (x^{\nu})^{\alpha_{\nu}}.$$ The
$j^{th}$
 vector of the standard basis is identified with $\partial_j$ the partial
derivative with respect to $x^j$. It will be convenient to denote
by $\hat{x}$ the sum of the $\nu$ variables $x^i$ and by
$\partial_{\hat{x}}$ the sum of the $\nu$ partial derivatives
$\partial_{x^i}$. Hence
$$\hat{x}=x^1+\cdots+x^{\nu}$$ and
$$\partial_{\hat{x}}=\partial_{x^1}+\cdots+\partial_{{x^\nu}}.$$
Observe that $\partial_{\hat{x}}\hat{x}=\nu$.

Let us denote by $V$ the vector space $\mathbb{R}^{\nu}$ and by
$S^{k}(V^{\star})$ its space of symmetric covariant tensors of
degree $k$. Any homogeneous polynomial vector field $V_k$ of
degree $k$ can be identified with an element of $V\bigotimes
S^{k}(V^{\star})$, and can be written as
$V_k=\sum_{i=1}^{\nu}V_k^i\partial_i$ where $V_k^i\in
S^{k}(V^{\star})$. The space $V\bigotimes S^{k}(V^{\star})$ admits
$$\{\frac{\mid\alpha\mid!}{\alpha !}x^{\alpha}\partial_i\}\mbox{ where }
i=1,\ldots,{\nu} \mbox{ and } \mid\alpha\mid=k$$ as basis. We
denote by $\parallel V_k\parallel_k$ the sup norm associated with
this basis, that is,
$$\parallel V_k\parallel_k=\max_{\alpha,i}\mid V^i_{k,\alpha}\mid$$
for any $V_k=\sum V^i_{k,\alpha}\frac{\mid\alpha\mid!}{\alpha
!}x^{\alpha}\partial_i$.

Let $\chi({\nu})=\chi_{-1}(\nu)$ denote the Lie algebra of formal
vector fields based at the origin $0$ of $\mathbb{R}^{\nu}$. For
every non-negative integer  $q$, let $\chi_{q}(\nu)$ denote the
Lie subalgebra of formal vector fields tangent up to order $q$ to
the zero vector field. This defines the decreasing filtration
naturally associated to $\chi(\nu)$. In component form any formal
vector field $V$ of $\chi_q(\nu)$ can be written as
$V=\sum_{i=1}^{\nu}V^i\partial_i$ with
$V^i=\sum_{\mid\alpha\mid>q}v^i_{\alpha}x^{\alpha}$. A formal
vector field $V$ will be said {\em positive} if its coefficients
$v^j_{\alpha}$ are positive.

Given a Lie pseudogroup $\Gamma^{\omega}$ acting in a neighborhood
of $0$ in $\mathbb{R}^{\nu}$, its formal Lie algebra ${\mathcal
L}(\Gamma)$ is a closed Lie subalgebra of $\chi({\nu})$ when the
latter is endowed with the Tychonov topology. Set ${\mathcal
L}_q(\Gamma)={\mathcal L}(\Gamma)\cap \chi_{q}({\nu})$. We have
${\mathcal L}(\Gamma)/{\mathcal L}_0(\Gamma)\simeq
\mathbb{R}^{\nu}$ whenever $\Gamma^{\omega}$ is transitive. We
associate to ${\mathcal L}(\Gamma)$ a  {\em flat} Lie algebra
denoted by $L(\Gamma)$. By definition, we have
$$L(\Gamma)=\oplus_{q=-1}^{\infty}{\mathcal L}_q(\Gamma)/{\mathcal
L}_{q+1}(\Gamma).$$

We note that the quotient space
$${\gamma}_q={\mathcal L}_{q-1}(\Gamma)/{\mathcal L}_{q}(\Gamma)$$
inherits from $V\bigotimes S^{k}(V^{\star})$ a Banach space
structure.

A transitive Lie pseudogroup is said to be {\em flat} if its Lie
algebra is isomorphic to its associated flat Lie algebra, i.e.
${\mathcal L}(\Gamma)\simeq L(\Gamma)$. It is well known that the
flat Lie algebra arising from a Lie pseudogroup $\Gamma^{\omega}$
characterizes the type of geometry associated to
$\Gamma^{\omega}$.  The pseudogroup of local transformations that
preserve a volume form and the Hamiltonian pseudogroup of local
symplectomorphisms of a symplectic manifold are classical examples
of flat Lie pseudogroups.

\subsection{Topology}\label{ntop}

Recall that formal vector field $V=\sum_{k=0}^{\infty}V_k$ (where,
for all $k$, $V_k\in V\bigotimes S^k(V^{\star})$) is analytic if
its coefficients satisfy $$\limsup
\parallel V_k\parallel_k^{\frac{1}{k}}<\infty.$$

Let $\rho$ be a positive real number and let ${\mathcal
L}_{\rho}^{\omega}(\Gamma)$ denote the subspace of ${\mathcal
L}^{\omega}(\Gamma)$ of $V$'s such that $\limsup \parallel
V_k\parallel_k/\rho^k<+\infty.$ We have ${\mathcal
L}^{\omega}(\Gamma)=\bigcup_{\rho>0}{\mathcal
L}_{\rho}^{\omega}(\Gamma)$. Each ${\mathcal
L}_{\rho}^{\omega}(\Gamma)$ is naturally endowed with a Banach
space structure with the norm
$$\parallel V\parallel_{\rho}=\sup_k  \frac{\parallel V_k\parallel_k}
{\rho^{k}}.$$ We shall mainly be concerned with isotropy
subgroups. In that case we will use preferably the equivalent norm
$$\parallel V\parallel_{\rho}=\sup_k  \frac{\parallel V_k\parallel_k}
{\rho^{k-1}}.$$

For $\rho<\rho'$ the injection ${\mathcal
L}_{\rho}^{\omega}(\Gamma)\hookrightarrow {\mathcal
L}_{\rho'}^{\omega}(\Gamma)$ is
 continuous and compact. Hence ${\mathcal L}^{\omega}(\Gamma)$ is a
complete Hausdorff locally convex topological
vector space. Its
 associated topology is the locally convex strict inductive limit topology
 $${\mathcal
L}^{\omega}(\Gamma)=\lim_{\stackrel{\longrightarrow}{\rho\in\mathbb{N}}}{\mathcal
L}_{\rho}^{\omega}(\Gamma).$$

 This endows ${\mathcal L}^{\omega}(\Gamma)$ with a Silva
topological Lie algebra structure. We shall denote by ${\mathcal
L}_{\rho,q}^{\omega}(\Gamma)$ the intersection
$${\mathcal L}_{\rho,q}^{\omega}(\Gamma)={\mathcal
L}_{\rho}^{\omega}(\Gamma)\cap\chi_q(\nu),$$ and by ${\mathcal
L}_{\rho,q,M}^{\omega}(\Gamma)$ the subset of ${\mathcal
L}_{\rho,q}^{\omega}(\Gamma)$ consisting of those analytic vector
fields $V$ satisfying
$$
\parallel V\parallel_{\rho}\leq M.
$$

\subsection{Statement of the main theorem}
Our main objective is to prove an infinite-dimensional version of
the classical second fundamental theorem of Lie for isotropy
subgroups of analytic Lie pseudogroups. Recall from \cite{KaR01},
Theorem 5.3, that the group $G_0^{\omega}(\nu)$ of analytic local
diffeomorphisms of $\mathbb{R}^{\nu}$ fixing the origin has the
natural structure of a G\^ateaux-analytic Lie group.

\begin{theo}[Lie II]\label{secondtheo} Any isotropy group of a Lie pseudogroup of
analytic transformations in $\nu$ variables is integrable into a
unique connected subgroup $H$ embedded in $G_0^{\omega}(\nu)$.
Such a subgroup $H$ is always a regular analytic Lie group
belonging to the class ${\mathcal CBH}^{\aleph_0}$ of
Campbell-Baker-Hausdorff groups of countable order.
\end{theo}

Let ${\mathcal L}^{\omega}(\Gamma)\subset \chi(\nu)$ be a Lie
algebra of local analytic vector fields defined in a neighborhood
of the origin $O$ of $\mathbb{R}^{\nu}$. For any integer $q$,
denote by ${\mathcal L}_q^{\omega}(\Gamma)$ the subalgebra of
vector fields in ${\mathcal L}^{\omega}(\Gamma)$ contained in
$\chi_q(\nu)$. Whenever $q'>q\geq 0$, the subalgebra ${\mathcal
L}_{q'}^{\omega}(\Gamma)$ is an ideal of ${\mathcal
L}_q^{\omega}(\Gamma)$. Suppose that for every integer $q\geq 1$,
we are given a linear section $\Sigma_q$ of the natural projection
$${\mathcal L}_{q-1}^{\omega}(\Gamma)\rightarrow{\mathcal
L}_{q-1}^{\omega}(\Gamma)/{\mathcal L}_q^{\omega}(\Gamma),$$ and
let
$$E_q=\Sigma_q({\mathcal L}_{q-1}^{\omega}(\Gamma)/{\mathcal
L}_q^{\omega}(\Gamma)).
$$

The isotropy Lie algebra ${\mathcal L}_0^{\omega}(\Gamma)$ then
decomposes as a direct sum of finite dimensional summands
$${\mathcal L}_0^{\omega}(\Gamma)=\bigoplus_{q=1}^{+\infty}E_q,$$
where In the sequel we will denote by $\Sigma$ the collection of
$\{\Sigma_q\}_{q\in\mathbb{N}}$. We shall prove the following {\em
strong boundedness property:}
\begin{lemm}\label{bound}
There exist linear sections $\Sigma_q:{\mathcal
L}_{q-1}^{\omega}(\Gamma)/{\mathcal
L}_q^{\omega}(\Gamma)\rightarrow{\mathcal
L}_{q-1}^{\omega}(\Gamma)$ taking values in the same Banach
subspace ${\mathcal L}_{\rho_0}^{\omega}(\Gamma)$ of ${\mathcal
L}^{\omega}(\Gamma)$, such that the norms of the continuous
operators $\Sigma_q:{\mathcal L}_{q-1}^{\omega}(\Gamma)/{\mathcal
L}_q^{\omega}(\Gamma)\rightarrow{\mathcal
L}_{q-1}^{\omega}(\Gamma)$ are uniformly bounded.
\end{lemm}
The preceding lemma will be used to prove the following theorem:

\begin{theo}[Lie III]\label{maintheo}
The isotropy Lie algebra ${\mathcal L}_0^{\omega}(\Gamma)$ is
integrable into a unique connected and simply connected analytic
Lie group $\Gamma_0$ of the second kind and countable order. In
other words $\Gamma_0$ is a regular analytic Lie group of class
${\mathcal CBH}^{\aleph_0}$.
\end{theo}

This latter theorem combined with Proposition~\ref{prop str boun}
of our Appendix implies Theorem~\ref{secondtheo}.

\section{Estimates}

We begin by deriving some basic estimates on the norm of the
exponential of a vector field $X\in {\mathcal
L}_{\rho,n}^{\omega}(\Gamma)$. These estimates will be needed to
prove the convergence of the infinite product of exponentials
which we will use to define the charts of Theorem~\ref{maintheo}.

We first consider $\Gamma^{\omega}$ the general Lie pseudogroup of
${\nu}$ variables acting on $\mathbb{R}^{\nu}$. It is clear that
$M\rho^n \frac{\hat{x}^{n+1}}{1-\rho \hat{x}}\partial_{\hat{x}}$
is the maximal positive vector field of ${\mathcal
L}_{\rho,n}^{\omega}(\Gamma)$ of norm less or equal to $M$. We
call such a vector field a canonical vector field with respect to
our bornological decomposition of the Lie algebra.

 Let $X=M\rho^n
\frac{\hat{x}^{n+1}}{1-\rho \hat{x}}\partial_{\hat{x}}$ and
$Y=N\rho^r\frac{\hat{x}^{r+1}}{1-\rho \hat{x}}\partial_{\hat{x}}$
be two canonical positive vector fields on $\mathbb{R}^{\nu}$.

Then
$$XY=MN\nu\rho^{n+r}\frac{\hat{x}^{n+r+1}}{1-\rho \hat{x}}\{
\frac{(r+1)}{1-\rho \hat{x}}+\frac{\rho \hat{x}}{(1-\rho
\hat{x})^2}\}\partial_{\hat{x}}$$ is also positive. Since in
addition
$$\frac{(r+1)}{1-\rho \hat{x}}+\frac{\rho \hat{x}}{(1-\rho
\hat{x})^2}=\frac{r+1-r\rho \hat{x}}
{(1-\rho \hat{x})^2},$$ we conclude that
$$XY\ll MN\nu\rho^{n+r}\hat{x}^{n+r+1}\frac{(r+1)}{(1-\rho
\hat{x})^3}\partial_{\hat{x}}.$$

More generally we obtain with
$Y=N\rho^r\frac{\hat{x}^{r+1}}{(1-\rho
\hat{x})^{\alpha}}\partial_{\hat{x}}$
$$XY=MN\nu\rho^{n+r}\frac{\hat{x}^{n+r+1}}{(1-\rho
\hat{x})^{\alpha+2}}\{r+1-(r+1-\alpha)\rho \hat{x}\}
\partial_{\hat{x}}$$
so that whenever $r+1\geq \alpha$
$$XY\ll MN\nu\rho^{n+r}\hat{x}^{n+r+1}\frac{(r+1)}{(1-\rho
\hat{x})^{\alpha+2}}\partial_{\hat{x}}.$$

We will deal with sequences of the form $X^nY=X(X^{n-1}Y)$ for $n$
integer. Let $X=M\rho^n \frac{\hat{x}^{n+1}}{1-\rho
\hat{x}}\partial_{\hat{x}}$ and
$Y=N\rho^r\frac{\hat{x}^{r+1}}{1-\rho \hat{x}}\partial_{\hat{x}}$
as before. Since
$$XY\ll MN\nu\rho^{n+r}\frac{\hat{x}^{n+r+1}}{(1-\rho
\hat{x})^3}(r+1)\partial_{\hat{x}}$$ we get
$$X^2Y=X(XY)\ll (\nu M)^2N\rho^{2n+r}\frac{\hat{x}^{2n+r+1}}{(1-\rho
\hat{x})^5}(r+1)(n+r+1)\partial_{\hat{x}}$$
whenever $n+r+1\geq 3$. This condition will be always satisfied
for $n\geq 2$. By recursion we get

\begin{lemm}
Let $X=M\rho^n \frac{\hat{x}^{n+1}}{1-\rho
\hat{x}}\partial_{\hat{x}}$ and
$Y=N\rho^r\frac{\hat{x}^{r+1}}{(1-\rho
\hat{x})^{\alpha}}\partial_{\hat{x}}$. For any integer $k$ and
$n\geq 2$
$$X^kY\ll (\nu M)^kN\rho^{kn+r}\frac{\hat{x}^{kn+r+1}}{(1-\rho
\hat{x})^{\alpha+2k}}(r+1)(n+r+1)
\cdots((k-1)n+r+1)\partial_{\hat{x}}.$$
\end{lemm}

Our next step is the study of $(id-X)\circ e^X$.

Remember that $$(id-X)\circ
e^X=id-\{\frac{1}{2!}X^2+\frac{2}{3!}X^3+\cdots+\frac{k}{(k+1)!}X^{k+1}+\cdots\}.$$
We consider $X=M\rho^n\frac{\hat{x}^{n+1}}{1-\rho
\hat{x}}\partial_{\hat{x}}$ with $n\geq 2$. Then
$$X^{k+1}\ll {\nu}^kM^{k+1}\rho^{(k+1)n}\frac{\hat{x}^{(k+1)n+1}}{(1-\rho
\hat{x})^{2k+1}}(n+1)(2n+1)
\cdots(kn+1)\partial_{\hat{x}}.$$ But $(n+1)(2n+1)\cdots(kn+1)\leq
(n+1)2(n+1)\cdots k(n+1)=k!(n+1)^k$ so that
$$\frac{1}{(k+1)!}X^{k+1}\ll
{\nu}^kM^{k+1}\rho^{(k+1)n}\frac{\hat{x}^{(k+1)n+1}}{(1-\rho
\hat{x})^{2k+1}}
\frac{(n+1)^k}{k+1}\partial_{\hat{x}}$$ and
$$\frac{k}{(k+1)!}X^{k+1}\ll \frac{k}{k+1}(\nu(n+1))^k
(M\rho^n)^{k+1}\frac{\hat{x}^{(k+1)n+1}}
{(1-\rho \hat{x})^{2k+1}}\partial_{\hat{x}}.$$

Notice now that all these vector fields belong to the Banach space
${\mathcal L}^{\omega}_{\sigma \rho}$ for any real number
$\sigma>1$. To evaluate their norm, we note that for any positive
integer $l$, we have
$$\frac{1}{(1-\rho y)^{l+1}}=\sum_{m=0}^{+\infty}\binom{m+l}{l} (\rho y)^m.$$

To obtain an upper bound of the norm of $\frac{k}{(k+1)!}X^{k+1}$
in ${\mathcal L}^{\omega}_{\sigma \rho}$, we start from
$$\frac{k}{(k+1)!}X^{k+1}\ll \frac{k}{k+1}(\nu(n+1))^k
M^{k+1}\frac{\rho^{n(k+1)}\hat{x}^{n(k+1)+1}}{(1-\rho
\hat{x})^{2k+1}}\partial_{\hat{x}}.$$ Put
$K=\frac{k}{k+1}(\nu(n+1))^k M^{k+1}$ and
$\rho_{\sigma}=\sigma\rho$ so $\rho=\sigma^{-1}\rho_{\sigma}$.
Therefore
$$\frac{k}{(k+1)!}X^{k+1}\ll
K\sigma^{-n(k+1)}\rho_{\sigma}^{n(k+1)}
\hat{x}^{n(k+1)+1}\left(\sum_{m=0}^{+\infty}\binom
{m+2k}{2k}\sigma^{-m}(\rho_{\sigma}\hat{x})^m\right)\partial_{\hat{x}}$$
It follows that the norm of $\frac{k}{(k+1)!}X^{k+1}$ in
${\mathcal L}^{\omega}_{\sigma \rho}$ is less than
$$K\sigma^{-n(k+1)}\max_{m\in\mathbb{N}} \{ \binom {m+2k}{2k}\sigma^{-m} \}.$$
Next, we have

\begin{lemm}\label{lmal}
For any positive real number $\alpha$, the sequence $
\{m^{\alpha}\sigma^{-m}\}_{m\in\mathbb{N}}$ admits
$\left(\frac{\alpha}{\ln\sigma}\right)^{\alpha}e^{-\alpha}$ as
upper bound.
\end{lemm}

\begin{proof}
Put $\phi(m)=m^{\alpha}\sigma^{-m}$. The derivative
$\frac{d\phi}{dm}=\sigma^{-m}m^{\alpha -1}(\alpha -m\ln\sigma)$
vanishes for $m_0=\frac{\alpha}{\ln\sigma}$. Since $m_0$
corresponds to a maximum, $\phi(m_0)$ provides an upper bound.
\end{proof}

Now remark that $\binom {m+2k}{2k}\leq \frac{(m+2k)^{2k}}{(2k)!}$,
so that
$$\binom {m+2k}{2k}\sigma^{-m}\leq \frac{(m+2k)^{2k}}{(2k)!}
\sigma^{-(m+2k)}\sigma^{2k}.$$
and
$$\max_{m\in\mathbb{N}}\{\binom {m+2k}{2k}\sigma^{-m}\}\leq
\max_{m\in\mathbb{N}}\{\sigma^{2k}\frac{(m+2k)^{2k}}{(2k)!}\sigma^{-(m+2k)}\}.$$
Using Lemma \ref{lmal} with $\alpha=2k$ we get
$$\max_{m\in\mathbb{N}}\{\binom{m+2k}{2k}\sigma^{-m}\}\leq\frac{\sigma^{2k}}{(2k)!}\left(\frac{2k}{\ln\sigma}\right)^{2k}e^{-2k}.$$
We obtain that way
$$\parallel \frac{k}{(k+1)!}X^{k+1} \parallel_{{\mathcal
L}^{\omega}_{\sigma\rho}}\leq
K\sigma^{-n(k+1)}\frac{\sigma^{2k}}{(2k)!}\left(\frac{2k}{\ln\sigma}\right)^{2k}e^{-2k}$$
that is

\begin{coro}
The norm of $\frac{k}{(k+1)!}X^{k+1}$ in ${\mathcal
L}^{\omega}_{\sigma, \rho}$ is bounded as follows
$$\parallel \frac{k}{(k+1)!}X^{k+1} \parallel_{{\mathcal
L}^{\omega}_{\sigma\rho}}\leq
M^{k+1}(\nu(n+1))^k\sigma^{-(n-2)k}e^{-2k}\left(\frac{2k}{\ln\sigma}\right)^{2k}\frac{\sigma^{-n}}{(2k)!}.$$
\end{coro}

Let us denote this upper bound by ${\mathcal B}(n,k,\sigma)$.
$${\mathcal
B}(n,k,\sigma)=\frac{M}{\sigma^n}\left(\sqrt{\frac{\nu(n+1)M}{\sigma^{n-2}}}\frac{2k}{e\ln\sigma}\right)^{2k}\frac{1}{(2k)!}.$$
From Stirling formula $n!>\frac{1}{2}(\frac{n}{e})^n\sqrt{2\pi n}$
so that $(2k)!>(\frac{2k}{e})^{2k}\sqrt{k\pi}$. Hence
$${\mathcal B}(n,k,\sigma)\leq
\frac{M}{\sigma^n}\left(\frac{\sqrt{\nu(n+1)M}}{{\sigma^{\frac{n-2}{2}}}\ln\sigma}\right)^{2k}\frac{1}{\sqrt{k\pi}}.$$

We conclude that the infinite sum
$$\frac{1}{2!}X^2+\frac{2}{3!}X^3+\cdots+\frac{k}{(k+1)!}X^{k+1}+\cdots$$
belongs to ${\mathcal L}^{\omega}_{\sigma, \rho}$ for any $\sigma$
satisfying ${\sigma^{\frac{n-2}{2}}}\ln\sigma>\sqrt{\nu(n+1)M}$.
Indeed if we put
$q=\frac{\sqrt{\nu(n+1)M}}{{\sigma^{\frac{n-2}{2}}}\ln\sigma}$,
$$\parallel\sum_{k=1}^{+\infty}\frac{k}{(k+1)!}X^{k+1}\parallel_{{\mathcal
L}_{\sigma\rho}}\leq
\frac{M}{\sigma^n}\sum_{k=1}^{+\infty}\frac{q^{2k}}{\sqrt{k\pi}}$$
which is less than
$\frac{M}{\sqrt{\pi}\sigma^n}\frac{q^2}{1-q^2}$.

We conclude this section on basic estimates by remarking that with
$X=M\rho^n\frac{\hat{x}^{n+1}}{1-\rho \hat{x}}\partial_{\hat{x}}$
and $n\geq 2$, we have
$$\frac{1}{k!}X^{k+1}\ll (\nu(n+1))^k
(M\rho^n)^{k+1}\frac{\hat{x}^{(k+1)n+1}}{(1-\rho
\hat{x})^{2k+1}}\partial_{\hat{x}}.$$ It follows that

\begin{coro}
For any $\sigma>1$ the norm of $\frac{1}{k!}X^{k+1}$ in ${\mathcal
L}^{\omega}_{\sigma ,  \rho}$ is bounded as follows
$$\parallel \frac{1}{k!}X^{k+1} \parallel_{{\mathcal L}^{\omega}_{\sigma \rho}}\leq
M^{k+1}(\nu(n+1))^k\sigma^{-(n-2)k}e^{-2k}\left(\frac{2k}{\ln\sigma}\right)^{2k}\frac{\sigma^{-n}}{(2k)!}.$$
\end{coro}

Remark that this upper bound is the same as the one given for
$\frac{k}{(k+1)!}X^{k+1}$ since we had replaced $\frac{k}{k+1}$ by
$1$. Therefore $X\circ e^X$ belongs to ${\mathcal L}_{\sigma\rho}$
whenever $q<1$ with
$q=\frac{\sqrt{\nu(n+1)M}}{{\sigma^{\frac{n-2}{2}}}\ln\sigma}$; we
also have
$$\parallel\sum_{k=1}^{+\infty}\frac{1}{k!}X^{k+1}\parallel_{{\mathcal
L}_{\sigma, \rho}}\leq
\frac{M}{\sqrt{\pi}\sigma^n}\frac{q^2}{1-q^2}.$$

\section{Iterative scheme and proof of Theorem 2}
We will assume in what follows that the uniform boundedness
condition of Lemma~\ref{bound} is satisfied. This assumption will
be justified in detail in the Appendix. It will be convenient to
represent the sections $\Sigma_q$ by the introduction of an
adapted filtered basis of the Lie algebra ${\mathcal
L}^{\omega}(\Gamma)$. A basis $\mathcal B$ of the Lie algebra
${\mathcal L}^{\omega}(\Gamma)$ will be said {\it filtered}
whenever, for all integer $q$, the subset ${\mathcal B}_q$ of its
elements that belong to ${\mathcal L}_q^{\omega}(\Gamma)$ forms a
basis of ${\mathcal L}_q^{\omega}(\Gamma)$. We have

\begin{lemm}
Any Lie algebra ${\mathcal L}^{\omega}(\Gamma)$ satisfying the
boundedness condition of Lemma~\ref{bound} admits a bounded
filtered basis $\mathcal B$. Such a basis can be generated by
lifting, via the sections $\Sigma_q$, any bounded homogeneous
basis of the associated flat Lie algebra $L^{\omega}(\Gamma)$.
More precisely there exists a constant $\rho_{0}>0$ such that for
any element $\varepsilon$ of ${\mathcal
L}_{q-1}^{\omega}(\Gamma)/{\mathcal L}_q^{\omega}(\Gamma)$, we
have $\parallel \tilde{\varepsilon}\parallel_{{\mathcal
L}_{\rho_{0}}}\leq
\parallel \varepsilon\parallel_q$ where
$\tilde{\varepsilon}=\Sigma_q(\varepsilon)$.
\end{lemm}

The quotient space ${\mathcal L}_q/{\mathcal L}_p$ for $p>q$ is
isomorphic
 to the
subspace of $\mathcal L$ spanned by the elements in ${\mathcal
B}_q$ that are not in ${\mathcal B}_p$. We will denote the image
of ${\mathcal L}^q/{\mathcal L}^p$ into ${\mathcal L}^q$ by the
section $\Sigma$ by $({\mathcal L}^q/{\mathcal L}^p)_{\mid
\Sigma}$.

  In order to obtain Theorem~\ref{maintheo} we need \cite{KaR01}
to demonstrate that any local transformation $\Phi$ of $\Gamma$ of
the form
$$\Phi=I-\sum_{i=1}^n(\sum_{\mid\alpha\mid\geq 2}\phi^i_{\alpha}
x^{\alpha})\partial_i,$$ with $I$ being the identity
transformation, can be written uniquely as an infinite product
$$
\cdots\circ\Exp\, v_n\circ\cdots\circ\Exp \,  v_2\circ\Exp \, v_1
$$
where  $v_n\in ({\mathcal L}^{\omega}_{2^n}(\Gamma)/{\mathcal
L}^{\omega}_{2^{n+1}}(\Gamma))_{\mid \Sigma}$ for each integer $n$
with the sum $\sum_n v_n$ being element of $\mathcal L$. The
existence and uniqueness of a formal solution
$(v_n)_{n\in\mathbb{N}}$ is easily established by iteration. The
core of the proof consists in checking the analyticity of the
series $\sum_n v_n$.

A basic ingredient in the procedure consists in the decomposition
of the ``free part" of any vector field with respect to the chosen
section $\Sigma$. We first illustrate the situation in the
one-dimensional case.

\subsection{Decompositions of analytic vector fields.}

\subsection{The one-dimensional case}

As a warm-up, we begin with the one-dimensional case by
considering analytic vector fields in a neighborhood of the origin
in $\mathbb{R}$.

Let first $p=2$ and consider a vector field $Z=(a_2x^3+a_3 x^4+a_4
x^5+\cdots)\partial_x$ in ${\mathcal L}^{\omega}_{\rho, 2,
M}(\Gamma)$. This means that its coefficients satisfy the
inequality
$$
\mid a_n\mid\leq M \rho^n.
$$
Suppose now that the adapted basis contains the vector fields
$$\epsilon_3=(b_2x^3+b_3x^4+b_4x^5+\cdots)\partial_x$$
and $$\epsilon_4=(c_3x^4+c_4x^5+\cdots)\partial_x.$$ Then we
decompose $Z$ as $Z=X+Y$ where $X$ represents the "free" part of
$Z$ and is given by
$$X=\frac{a_2}{b_2}\epsilon_3+(a_3-\frac{a_2}{b_2}b_3)\frac{1}{c_3}\epsilon_4$$
and where $Y$ is given by
$Y=(\alpha_4x^5+\alpha_5x^6+\cdots+\alpha_nx^{n+1}+\cdots)\partial_x$
with
$$\alpha_n=a_n-a_2\frac{b_n}{b_2}-a_3\frac{c_n}{c_3}+a_2\frac{b_3c_n}{b_2c_3}.$$

Let now $p=3$ and consider a vector field $Z=(a_3 x^4+a_4
x^5+a_5x^6+a_6x^7\cdots)\partial_x$ in ${\mathcal
L}^{\omega}_{\rho,3,M}(\Gamma)$. Suppose that the adapted basis
contains the vector fields
$$\epsilon_4=(c_3x^4+c_4x^5+c_5x^6+c_6x^7\cdots)\partial_x,$$
$$\epsilon_5=(d_4x^5+d_5x^6+d_6x^7\cdots)\partial_x,$$
$$\epsilon_6=(f_5x^6+f_6x^7\cdots)\partial_x.$$
We now have a decomposition $Z=X+Y$ where the "free" part $X$ of
$Z$ takes the form
$$X=\frac{a_3}{c_3}\epsilon_4+(a_4-\frac{a_3}{c_3}c_4)\frac{1}{d_4}\epsilon_5+
\left(a_5-\frac{a_3}{c_3}c_5-(a_4-\frac{a_3}{c_3}c_4)\frac{d_5}{d_4}\right)\frac{1}{f_5}\epsilon_6$$
and where
$Y=(\alpha_6x^7+\alpha_7x^8+\cdots+\alpha_nx^{n+1}+\cdots)\partial_{x}$
with
$$\alpha_n=a_n-\frac{a_3}{c_3}c_n-(a_4-\frac{a_3}{c_3}c_4)\frac{d_n}{d_4}-\left(a_5-\frac{a_3}{c_3}c_5-(a_4-\frac{a_3}{c_3}c_4)\frac{d_5}{d_4}\right)\frac{f_n}{f_5}.$$

By our strong boundedness hypothesis we can find a basis
$\{\epsilon_i\}$ for which $\frac{b_n}{b_2}\leq {\rho_0}^{n-2}$,
$\frac{c_n}{c_3}\leq {\rho_0}^{n-3}$ etc. for some positive real
number $\rho_0>0$. It is then always possible to reduce the
problem to the case where $\rho_0=1$ using the adjoint action
$Ad(\rho_0x)$.

When $p=2$, we obtain $\mid\alpha_n\mid\leq \mid a_n\mid+\mid
a_3\mid+2\mid a_2\mid$ and when $p=3$, we have
$\mid\alpha_n\mid\leq\mid a_n\mid+\mid a_5\mid+2\mid a_4\mid+4\mid
a_3\mid$. In the general case we will get
$$\mid\alpha_n\mid\leq\mid a_n\mid+\mid a_{2p-1}\mid+2\mid
a_{2p-2}\mid+2^2\mid a_{2p-3}\mid+\cdots+2^{p-1}\mid a_p\mid.$$

Now $\mid a_n\mid\leq M \rho^n$ so that
$$\mid\alpha_n\mid\leq
M\left(\rho^n+\rho^{2p-1}(1+2x+2^2x^2+\cdots+2^{p-1}x^{p-1})\right)$$
with $x=\frac{1}{\rho}$. But
$1+2x+2^2x^2+\cdots+2^{p-1}x^{p-1}+\cdots=\frac{1}{1-2x}$ for
$2x<1$. Hence if we start at $\rho\geq 4$ each $\mid\alpha_n\mid$
is bounded over by $M(\rho^n+2\rho^{2p-1})$ so that
$$\parallel Y\parallel_{{\mathcal L}_{\rho}}\leq M \limsup_{k\geq
1}\{1+\frac{2}{\rho^k}\}=\frac{3}{2}M.$$
Since $X=Z-Y$ we have proved the following

\begin{prop}
With a basis $\{\epsilon_i\}$ for which $\frac{b_n}{b_2}\leq 1$,
$\frac{c_n}{c_3}\leq 1, \ldots$ for all $n$ the decomposition
splits any $Z\in{\mathcal L}_{\rho,p,M}(\Gamma)$ into $X+Y$ where
$X\in{\mathcal L}^{\omega}_{\rho,p,\frac{5}{2}M}(\Gamma)$ and
$Y\in{\mathcal L}^{\omega}_{\rho, 2p, \frac{3}{2}M}(\Gamma)$.
\end{prop}

\subsection{The higher dimensional case.}

  Let $Z=Z_{p+1}+Z_{p+2}+\cdots$ be an element of ${\mathcal
L}^{\omega}_{\rho,p,M}(\Gamma)$. Here each $Z_{p+i}$ represents
the homogeneous part of $Z$ of degree $p+i$. For each integer $i$
denote by $\tilde{Z}_{p+i}$ the image of $Z_{p+i}$ by the chosen
section $\Sigma$; i.e.
$\tilde{Z}_{p+i}=\Sigma(Z_{p+i})=\Sigma_{p+i}(Z_{p+i})$. We can
rewrite $Z$ as follows
$$Z=\tilde{Z}_{p+1}+(Z_{p+2}-\tilde{Z}^{p+2}_{p+1})+
(Z_{p+3}-\tilde{Z}^{p+3}_{p+1})+\cdots$$ where, for each integer
$l>i$, $\tilde{Z}^{p+l}_{p+i}$ is the homogeneous part of degree
$p+l$ of $\tilde{Z}_{p+i}$.

Since $Z\in {\mathcal L}^{\omega}_{\rho,p,M}(\Gamma)$, its
homogeneous part of degree $p+1$ satisfies $\parallel
Z_{p+1}\parallel_{p+1}\leq M\rho^p$. Therefore
$$\parallel \tilde{Z}^{p+l}_{p+1}\parallel_{p+l}\leq  M\rho^p
\rho_0^{l-1}\leq M\rho^p$$ for all integers $l\geq 1$ whenever
$\rho_0$ has been chosen smaller than $1$ with no loss of
generality. Denote now by $\tilde{Z}_{p+2}$ the image under
$\Sigma$ of the homogeneous part $Z_{p+2}-\tilde{Z}^{p+2}_{p+1}$
of degree $p+2$.

With this notation $Z$ takes the form
$$Z=\tilde{Z}_{p+1}+\tilde{Z}_{p+2}+
\sum_{l=3}^{+\infty}(Z_{p+l}-\tilde{Z}^{p+l}_{p+1}-\tilde{Z}^{p+l}_{p+2}).$$

At that stage we have the following bounds
$$\parallel \tilde{Z}^{p+2}_{p+2}\parallel_{p+2}=
\parallel Z_{p+2}-\tilde{Z}^{p+2}_{p+1}\parallel_{p+2}\leq
\parallel Z_{p+2}\parallel_{p+2}+\parallel
\tilde{Z}^{p+2}_{p+1}\parallel_{p+2},$$
that is
$$\parallel \tilde{Z}^{p+2}_{p+2}\parallel_{p+2}\leq
M\rho^{p+1}+ M\rho^p,$$ so that $$\parallel
\tilde{Z}^{p+l}_{p+2}\parallel_{p+2}\leq M\rho^{p+1}+ M\rho^p,$$
for all integers $l\geq 2$.

Continuing this way we obtain the decomposition of $Z=X+Y$ into a
free part $X$ and a remainder $Y$. The free part $X$ takes the
form
$$X=\tilde{Z}_{p+1}+\tilde{Z}_{p+2}+\cdots+
\tilde{Z}_{2p}$$ while the remainder $Y$ is
$$Y=\sum_{l=p+1}^{+\infty}(Z_{p+l}-\tilde{Z}^{p+l}_{p+1}-\tilde{Z}^{p+l}_{p+2}-\cdots-\tilde{Z}^{p+l}_{2p})$$

Moreover
$$\parallel Y_n\parallel_n\leq \parallel Z_n\parallel_n+\parallel
\tilde{Z}^{n}_{p+1}\parallel_n+
\parallel \tilde{Z}^{n}_{p+2}\parallel_n+\cdots+\parallel
\tilde{Z}^{n}_{2p}\parallel_n$$
so that
$$\parallel Y_{n+1}\parallel_{n+1}\leq M\rho^n +M\rho^p
+(M\rho^{p+1}+M\rho^p)+\cdots+(M\rho^{2p-1}+\cdots+M\rho^p)$$
or
$$\parallel Y_{n+1}\parallel_{n+1}\leq
M\left(\rho^n+\rho^{2p-1}(1+2x+2^2x^2+\cdots+2^{p-1}x^{p-1})\right)$$
with $x=\frac{1}{\rho}$. But
$1+2x+2^2x^2+\cdots+2^{p-1}x^{p-1}+\cdots=\frac{1}{1-2x}$ for
$2x<1$. Hence if we start at $\rho\geq 4$ each $\parallel
Y_{n+1}\parallel_{n+1}$ is bounded over by
$M(\rho^n+2\rho^{2p-1})$. The norm of the remainder $Y$ in
${\mathcal L}_{\rho}$ is bounded by
$$\parallel Y\parallel_{{\mathcal L}_{\rho}}\leq M \limsup_{k\geq
1}\{1+\frac{2}{\rho^k}\}\leq\frac{3}{2}M.$$

Since $X=Z-Y$ we have proved the following

\begin{prop}
The decomposition splits any $Z\in{\mathcal
L}^{\omega}_{\rho,p,M}(\Gamma)$ into $X+Y$ where $X\in{\mathcal
L}^{\omega}_{\rho,p,\frac{5}{2}M}(\Gamma)$ and $Y\in{\mathcal
L}^{\omega}_{\rho,2p,\frac{3}{2}M}(\Gamma)$.
\end{prop}

\subsection{Proof of Theorem 2 - General iterative scheme}

The iterative step makes us pass from $x-Z$ to $(x-X)\circ
e^X-Y\circ e^X=x-W$ where $Z=X+Y$ is the adapted decomposition
described above. If $Z$ belongs to ${\mathcal
L}^{\omega}_{\rho,p,M}(\Gamma)$ then $X$ belongs to ${\mathcal
L}^{\omega}_{\rho,p, M'=\frac{5}{2}M}(\Gamma)$ and $Y$ to
${\mathcal L}^{\omega}_{\rho,2p, M''=\frac{3}{2}M}(\Gamma)$.
Therefore we obtain
$$\parallel W\parallel_{{\mathcal L}^{\omega}_{\sigma\rho}}\leq
\frac{2M'}{\sqrt{\pi}\sigma^p}\frac{q^2}{1-q^2}+\parallel
Y\parallel_{{\mathcal L}^{\omega}_{\sigma\rho}},$$ where
$q=\frac{\sqrt{(p+1)M'}}{\sigma^{\frac{p-2}{2}}\ln\sigma}$ has
been chosen  to be less than $1$. Since $Y$ is of order $2p$, that
is $Y\in {\mathcal L}^{\omega}_{\rho,2p,M''}(\Gamma)$, we have
also $\parallel Y\parallel_{{\mathcal
L}^{\omega}_{\sigma\rho}}\leq \frac{M''}{\sigma^{2p}}.$ Finally
$$\parallel W\parallel_{{\mathcal L}^{\omega}_{\sigma\rho}}\leq
\frac{5M}{\sqrt{\pi}\sigma^p}\frac{q^2}{1-q^2}+\frac{3M}{2\sigma^{2p}}.$$

The iteration sends ${\mathcal L}^{\omega}_{\rho,p,M}(\Gamma)$ to
${\mathcal L}^{\omega}_{\sigma\rho,2p,\tilde{M}}(\Gamma)$ where
$$q=\frac{\sqrt{(p+1)5M/2}}{\sigma^{\frac{p-2}{2}}\ln\sigma}<1$$ and
$$\tilde{M}=\frac{M}{\sigma^p}\left(\frac{5}{\sqrt{\pi}}\frac{q^2}{1-q^2}+\frac{3}{2\sigma^p}\right).$$

With no loss of generality, we can set the value $q$ to $1/2$. We
can use the "LambertW" function solve the equation
$$\sigma^{\frac{p-2}{2}}\ln\sigma=\sqrt{10(p+1)M}$$
for $\sigma$ as a function of $p$ and $M$. This function is
defined as the principal branch of the solution of the equation
$z=w\,exp(w)$ for $w$ as a function of $z$. The Lambert function
$\LW(z)$ satisfies the differential equation
$$dw/dz=w/(z(1+w)).$$ More generally, the solution of $x^a\ln(x)=b$ is
given by $x=e^{\frac{\LW(ab)}{a}}$, and we obtain
$$\sigma=e^{2\frac{\LW(\frac{p-2}{2}\sqrt{10(p+1)M})}{(p-2)}}.$$

The iteration takes us from ${\mathcal
L}^{\omega}_{\rho,n_{i},M_{i}}(\Gamma)$ to ${\mathcal
L}^{\omega}_{\rho_{i+1},n_{i+1},M_{i+1}}(\Gamma)$ where
$n_{i+1}=2n_i$, $\rho_{i+1}=\sigma_i\rho_i$,
$$M_{i+1}=\frac{M_i}{\sigma_i^{n_i}}\left(\frac{5}{3\sqrt{\pi}}+\frac{3}{2{\sigma_i}^{n_i}}\right)$$
and
$$\sigma_i=e^{2\frac{\LW(({n_i}/2-1)\sqrt{10(n_i+1)M_i})}{(n_i-2)}}.$$

Obviously $n_i=2^in_0$. Moreover since
$\frac{5}{3\sqrt{\pi}}\approx 0.94$ and $\sigma_i>1$ for all
integers $i$ it is clear that $M_i\leq (\frac{5}{2})^iM_0$.
Therefore
$$\rho_{i+1}=\sigma_i\sigma_{i-1}\cdots\sigma_0\rho_0$$
is less than
$$\rho_{i+1}\leq
e^{\sum_{k=0}^i2\frac{\LW((2^{k-1}n_0-1)\sqrt{10(2^kn_0+1)(5/2)^kM_0})}{(2^kn_0-2)}}\rho_0$$
But according to the asymptotic properties of $\LW$, which are
similar to those of the
 classical logarithmic function, the infinite sum
$$\sum_{k=0}^{+\infty}2\frac{\LW((2^{k-1}n_0-1)\sqrt{10(2^kn_0+1)(5/2)^kM_0})}{(2^kn_0-2)}$$
converges. This allows to prove the bornological convergence of
the iteration completing the proof of Theorem~\ref{maintheo}.

\section{Appendix: The Malgrange estimate}
In this appendix we prove that the Lie algebra of a Lie
pseudogroup of analytic local transformations satisfies the
boundedness condition of Lemma~\ref{bound}. This property is
obtained using Malgrange's {\it Appendice sur le th\'eor\`eme de
Cartan-K\"ahler} \cite{Mal72} [part II p.134] to which the reader
is referred.

\subsection{General facts}
\subsubsection{Notation}
Throughout this section we will consider a system of $r$ equations
with $p$ independent variables and $q$ dependent ones. The
independent variables will be respectively denoted by
$z=(z^1,\ldots,z^p)$, and the dependent ones by
$y=(y^1,\ldots,y^q)^T$. Note that $y$ will be viewed as a column
vector. For each $i\in\{1,\ldots,p\}$ and $j\in\{1,\ldots,q\}$ the
partial derivative $\partial_{z^i}y^j$ will be denoted by $y^j_i$.
We put $y_i=(y_i^1,\ldots,y_i^q)^T$. We will use small cap greek
letters to represent multi-indices of format $p$. If
$\alpha=(\alpha_1,\ldots,\alpha_p)$ is such a multi-index, its
entries $\alpha_i$ are nonnegative integers. The greek letter
$\varepsilon$ will be used specifically as follows; we will denote
by
$$\varepsilon_i=(0,\ldots,0,1,0,\ldots,0)$$ the $p$-multiplet
with
$1$ in $i$-th position and $0$ elsewhere. The null $p$-multiplet ($0$
everywhere) will be represented by $\bf 0$. As usual
$\mid\alpha\mid=\alpha_1+\cdots+\alpha_p$ (length of $\alpha$),
$\alpha!=\alpha_1!\cdots\alpha_p!$ and
$z^{\alpha}=(z^1)^{\alpha_1}\cdots(z^p)^{\alpha_p}$.

\subsubsection{Defining equations}
A general linear and homogeneous system of first-order PDEs, that
we denote by $(E)$, can be written in matrix form as
$$A(z)y+\sum_{i=1}^{p}B^i(z)y_i=0$$
where the $r\times q$-matrices $A(z), B^i(z)$ are analytic in $z$.
We are working in a neighborhood of the origin for the independent
variable $z$, $A(z), B^i(z)$ will be represented by their Taylor
series centered at the origin i.e.
$$A(z)=\sum_{\mid\alpha\mid\geq 0} A_{\alpha}z^{\alpha}$$
and
$$B^i(z)=\sum_{\mid\alpha\mid\geq 0}B^i_{\alpha}z^{\alpha}$$
where $A_{\alpha}, B^i_{\alpha}$ are $r\times q$-matrices with constant
coefficients (real or complex).

From now on we will use the summation convention.
\subsection{Recursion}
\subsubsection{Formal solution}
We will deal with formal solutions $y$ of the form
$F=f_{\gamma}z^{\gamma}$. Such an $F$ is constructed step-by-step
by replacing $y$ by $F$ in the equation and solving the
corresponding equation at order $0$, then $1$, then $2$ and so on.
We focus only on strongly extendible jets ( {\it jets fortements
prolongeables} in Malgrange's terminology \cite{Mal72}[part II
p.138] ). The recursive construction is therefore made possible as
soon as $f_{\gamma}$ for $\mid\gamma\mid=0,1$ have been chosen by
solving the equation $(E)$ at order $0$, i.e. by solving what we
will call $(E_0)$.

\subsubsection{Equation at order $l$}
From $y=f_{\gamma}z^{\gamma}$ we
get $y_i=(\gamma_i+1)f_{\gamma+\varepsilon_i}z^{\gamma}$ where $\gamma_i$ is
the $i$-th entry of $\gamma=(\gamma_1,\ldots,\gamma_p)$. The equation $(E)$
at order $l$, that we denote by $(E_l)$, takes the matrix form
$$\sum_{{\mid\beta\mid=l}}B^i_{\bf 0
}(\beta_i+1)f_{\beta+\varepsilon_i}z^{\beta}=
-[\sum_{\mid\beta\mid=l}A_{\alpha}f_{\beta-\alpha}z^{\beta}+
\sum_{\stackrel{\mid\beta\mid=l}{\beta\neq\gamma}}B^i_{\beta-\gamma}(\gamma_i+1)f_{\gamma+\varepsilon_i}z^{\beta}]$$
where the summation over $i$ is implicit. Equation $(E_l)$ is
precisely, in our case, equation
$$\sum_{\stackrel{\mid\alpha\mid=k}{\mid\beta\mid=l}}u_{\alpha}f_{\alpha+\beta}\frac{(\alpha+\beta)!}{\beta!}z^{\beta}=\Psi(z,D^{\alpha}F^{k+l-1})_l$$
of Malgrange \cite{Mal72}[part II p.139]. We will denote this
equation by $(\Phi_l)$ and, following Malgrange, we will denote by
$G_l=\sum_{\mid\beta\mid=l} g_{\beta} z^{\beta}$ its right-hand
side.

\subsection{Norms and the Malgrange estimate}

\subsubsection{Norms}
In accordance with Malgrange's notation we will denote by $E_1$
the space $\mathbb{R}^r$ (or $\mathbb{C}^r$) and by $E$ the space
$\mathbb{R}^q$ (or $\mathbb{C}^q$). The space of $r\times
q$-matrices (over $\mathbb{R}$ or $\mathbb{C}$) will simply by
denoted by $\mathcal M$. We endow these three vector spaces with
the sup norm. Thus for
$f_{\alpha}=(f_{\alpha}^1,\ldots,f_{\alpha}^q)^T$ in $E$, $$\mid
f_{\alpha}\mid= \max_{i=1,\ldots,q}\mid f_{\alpha}^i\mid$$ and
likewise, for $A_{\alpha}=[(A_{\alpha})^i_j]$ in $\mathcal M$,
$$\mid A_{\alpha}\mid=\max_{\stackrel{i=1,\ldots,r}{j=1,\ldots,q}}\mid
(A_{\alpha})^i_j\mid.$$

In addition we will endow the space $S^l\bigotimes E$ (resp. $S^l\bigotimes
E_1$ and $S^l\bigotimes {\mathcal M}$) of degree $l$ homogeneous polynomials
(in the $z$ variable) with values in $E$ (resp. $E_1$ and $\mathcal M$) with
the following norm; if $H=\sum_{\mid\beta\mid=l}h_{\beta}z^{\beta}$ is such
an element then we put
$$\parallel H\parallel_l=\max_{\mid\beta\mid=l}\frac{\beta!}{l!}\mid
h_{\beta}\mid$$ where $\mid
h_{\beta}\mid$ is the corresponding norm of $h_{\beta}$ in $E$ (resp. $E_1$
and $\mathcal M$).

The following proposition shows to which extent these norms are well adapted.
\begin{prop}\label{multiplicative norm}
If $A=A_l(z)$ is a homogeneous polynomial of degree $l$ with values in
$\mathcal M$ and if $Y=Y_m(z)$ is a homogeneous polynomial of degree $m$
with values in $E$ (column matrix of format $q$) then $AY$ is a homogeneous
polynomial of degree $l+m$ with values in $E_1$ (column matrix of format
$r$) satisfying $$\parallel AY\parallel_{l+m}\leq q\parallel
A\parallel_l\parallel Y\parallel_m$$
\end{prop}

The proof of the preceding proposition is an adaptation of
Proposition 3.3 \cite{Mal72}[part II p.137] taking into account
our choice of norms.

\subsubsection{The Malgrange estimate}
With our choice of norms, the Malgrange estimate \cite{Mal72}[part
II p.139],
$$\max_{\mid \alpha\mid=k+l}\alpha!\mid f_{\alpha}\mid r^{\alpha}\leq C
\max_{\mid\beta\mid=l}\beta!\mid g_{\beta}\mid r^{\beta}$$
obtained from equation $(\Phi_l)$ as a consequence of his
Corollary 2.2 \cite{Mal72}[part II p.135] reads
$$(k+l)!\parallel F\parallel_{k+l}\leq \frac{C}{r_0^k}l!\parallel
G_l\parallel_l$$ if we consider as a system of $p$ strictly
positive numbers\footnote{Note that the integer $n$ in Malgrange
paper coincides with our $p$} $r=(r_1,\ldots,r_p)$ the string
$r=(r_0,\ldots,r_0)$ ( see Section 2 on {\it Polydisques
distingu\'es} in \cite{Mal72}[part II p.134] ).

In our particular case of first order PDE, this inequality becomes
$$\parallel F\parallel_{l+1}\leq\frac{C}{r_0(l+1)}\parallel G_l\parallel_l$$
which we denote by $({\mathcal I}_{l})$.
\begin{rema}
$F=f_{\alpha}z^{\alpha}$ is analytic if and only if there exist $M,\rho>0$
such that $\parallel F\parallel_l\leq M\rho^l$ for all integer $l$. This is
a consequence of the inequality
$$\mid\sum_{\mid\alpha\mid=l}f_{\alpha}z^{\alpha}\mid\leq\sum_{\mid\alpha\mid=l}\mid
f_{\alpha}\mid\mid z\mid^{\alpha}\leq
(\sum_{\mid\alpha\mid=l}\frac{l!}{\alpha !}\mid z\mid^{\alpha})\parallel
F\parallel_l\leq \parallel F\parallel_l (\mid z^1\mid+\cdots+\mid
z^p\mid)^l.$$
\end{rema}

\subsection{Prolongation - Strong boundedness}
\subsubsection{Growth of coefficients}
From Proposition~\ref{multiplicative norm} the norm of the
right-hand side $G_l=\sum_{\mid\beta\mid=l} g_{\beta} z^{\beta}$
of $(E_l)$ satisfies the inequality
$$\parallel G_l\parallel_l\leq q[\sum_{s+t=l}\parallel A\parallel_s\parallel
F\parallel_t+\sum_{\stackrel{s+t=l}{s>0}}\parallel B^i\parallel_s\parallel
F_i\parallel_t]$$
in which intervene the homogeneous parts
$(F_i)_t=\sum_{\mid\gamma\mid=t}(\gamma_i+1)f_{\gamma+\varepsilon_i}
z^{\gamma}$
of $F_i$. Since the length $\mid\gamma+\varepsilon_i\mid=t+1$ for
$\mid\gamma\mid=t$ and
$$\max_{\mid\alpha\mid=t+1}\frac{\alpha!}{\mid\alpha\mid!}\mid
f_{\alpha}\mid=\parallel F\parallel_{t+1}$$ we obtain, for
$\mid\alpha\mid=t+1$,
$$\mid f_{\alpha}\mid\leq\frac{(t+1)!}{\alpha!}\parallel F\parallel_{t+1}.$$
Therefore $$\mid
f_{\gamma+\varepsilon_i}\mid\leq\frac{(t+1)!}{(\gamma_i+1)\gamma!}\parallel
F\parallel_{t+1}$$ so that
$$\frac{\gamma!}{t!}\mid(\gamma_i+1)f_{\gamma+\varepsilon_i}\mid\leq
(t+1)\parallel F\parallel_{t+1}$$ for all multi-indices $\gamma$
of length $t$. As a consequence $$\parallel F_i\parallel_t\leq
(t+1)\parallel F\parallel_{t+1}$$ and finally, introducing two
constants $M,\rho_0>0$ such that $\parallel A\parallel_s,\parallel
B^i\parallel_s\leq M\rho_0^s$ for all integer $s$, we get
$$\parallel G_l\parallel_l\leq q[\sum_{s+t=l}M\rho_0^s\parallel
F\parallel_t+\sum_{\stackrel{s+t=l}{s>0}}M\rho_0^s(t+1)\parallel
F\parallel_{t+1}].$$ We deduce from $({\mathcal I}_l)$ the
following lemma

\begin{lemm}\label{lemm gofcoef} The prolongation scheme
described by Malgrange leads, in the case of first order linear
and homogeneous PDEs, to the inequality
$$\parallel F\parallel_{l+1}\leq \frac{CMq}{r_0(l+1)}[(\rho_0^l\parallel
F\parallel_0+\rho_0^{l-1}\parallel F\parallel_1+\cdots+\parallel
F\parallel_l)$$$$+(\rho_0 l\parallel F\parallel_l+\rho_0^2(l-1)\parallel
F\parallel_{l-1}+\cdots+\rho_0^l\parallel F\parallel_0)]
.$$
\end{lemm}

\subsubsection{Strong boundedness}

\begin{prop}[Strong boundedness]\label{prop str boun}
Let $\Lambda$ be a linear and homogeneous analytic system of
partial differential equations. Then its vector space of analytic
local solutions admits as basis a {\it bounded} family $\mathcal
B$ satisfying the following growth condition: for all $F$ in
$\mathcal B$ and all positive integers $k$
$$\parallel F\parallel_{d+k}\leq K^k\parallel F\parallel_d$$
where $K>0$ is a constant depending only on $\Lambda$,  and $d$
represents the smallest degree appearing in $F$.
\end{prop}
\begin{proof}
Such a basis will be constructed using Malgrange prolongation
scheme. For instance the ``affine part" will be obtained by
prolongation of a basis of the vector space solution to $(E_0)$.
The ``homogeneous part" of degree 2 will be obtained by
prolongation of a basis of the vector space solution to $(E_1)$
without right-hand side (no affine part) and so on.

In order to control the growth of the coefficients we associate to
each such prolongation $F$ ( element of the basis $\mathcal B$ ) a
sequence $(\phi_n)$ as follows. For a prolongation representing
the affine part we put $\phi_0=\parallel F\parallel_0$,
$\phi_1=\parallel F\parallel_1$ and for a prolongation we put
$F=F_d+F_{d+1}+\cdots$ representing the homogeneous part of degree
$d>1$  put $\phi_l=0$ for $l<d$ and $\phi_d=\parallel
F_d\parallel_d=\parallel F\parallel_d$. Finally for all integers
$l\geq 1$ set
$$\phi_{l+1}=\frac{CMq}{r_0(l+1)}[(\rho_0^l\phi_0+\rho_0^{l-1}\phi_1+
\cdots+\phi_l)+(\rho_0l\phi_l+\rho_0^2(l-1)\phi_{l-1}\cdots\rho_0^l\phi_0)].$$
The sequence $(\phi_n)_{n\in\mathbb{N}}$ bounds the sequence
$(\parallel F\parallel_n)_{n\in\mathbb{N}}$ from above, i.e.
$\parallel F\parallel_n\leq\phi_n$ for all integer $n$. Moreover
$$\phi_{l+1}=\frac{CMq}{r_0(l+1)}[\rho_0(\rho_0^{l-1}\phi_0+
\cdots+\phi_{l-1})+\phi_l+\rho_0l\phi_l+\rho_0(\rho_0(l-1)\phi_{l-1}\cdots
\rho_0^{l-1}\phi_0)],$$
that is to say
$$\phi_{l+1}=\frac{CMq}{r_0(l+1)}[\rho_0(\frac{lr_0\phi_l}{CMq})+\phi_l+
\rho_0l\phi_l],$$
or equivalently
$$\phi_{l+1}=[\frac{l}{l+1}\rho_0+\frac{CMq}{r_0(l+1)}(1+\rho_0l)]\phi_l.$$
If we put $K=\rho_0+\frac{CMq}{r_0}+\frac{CMq}{r_0}\rho_0$ this gives rise to
$$\phi_{l+1}\leq K\phi_l$$ for all integer $l\geq 1$.

Therefore $$\parallel F\parallel_{d+k}\leq\phi_{d+k}\leq
K^k\phi_d=K^k\parallel F\parallel_d$$ as claimed. This inequality
holds at least for $d\geq 2$. We include the cases $d=0$ and $d=1$
by replacing if necessary $K$ by a larger constant.
\end{proof}


\begin{thebibliography}{2}%{9}

\bibitem[Ca 04]{Ca04} E.~Cartan, {\em Sur la structure des groupes infinis de
transformations}, Ann. Ec. Normale \textbf{21} (1904), 219-308.

\bibitem[Ca 37]{Ca37} E.~Cartan, {\em La structure des groupes
infinis}, Seminaire de Math., expos\'e G, 1er mars 1937, reprinted
in {\em Elie Cartan, Oeuvres compl\`etes, Vol. II}, Editions du
CNRS, 1984.

\bibitem[Ch 53]{Ch53} S.S. ~Chern, {\em Pseudogroupes continus infinis},
Colloque de g\'eom\'etrie diff\'erentielle de Strasbourg, Editions
du CNRS, 1953, reprinted in {\em Shiing-Shen Chern, Selected
papers}, Springer-Verlag, 1978.

\bibitem[GuS 64]{guillemin 64} V.~Guillemin, S.~Sternberg, {\em An algebraic model of transitive differential geometry,} Bull.
of the A.M.S., vol. \textbf{70}, No. 1 (1964), 16--47.

\bibitem[KR 97-1]{KaR97-1} N.~Kamran, T.~Robart, {\em Abstract structure for
Lie pseudogroups of infinite type}, C. R. Acad. Sci. Paris, t.
\textbf{324}, S\'erie I (1997), 1395--1399.

\bibitem[KR 97-2]{KaR97-2}N.~Kamran, T.~Robart, {\em Perspectives sur la th\'eorie des pseudogroupes de transformations
analytiques de type infini}, J. Geom. Phys. \textbf{23} (1997),
no. 3-4, 308--318.

\bibitem[KR 00]{KaR00} N.~Kamran, T.~Robart, {\em On the parametrization problem of Lie pseudogroups of infinite
type}, C. R. Acad. Sci. Paris Sér. I Math. \textbf{331} (2000),
no. 11, 899--903.

\bibitem[KR 01]{KaR01} N.~Kamran, T.~Robart, {\em A manifold structure for analytic
isotropy {Lie}  pseudogroups of infinite type}, Journal of Lie
Theory \textbf{11} (2001) no. 1, 57--80.

\bibitem[KS 72]{KS72}A.~Kumpera, D.~Spencer, {\em Lie equations. Vol. I: General theory}, Annals of Mathematics Studies,
\textbf{73}, Princeton University Press, 1972.

\bibitem[L 82]{L82}J.~Leslie, {\em On the group of real analytic diffeomorphisms of a compact real analytic manifold},
Trans. Amer. Math. Soc. \textbf{274} (1982), no. 2, 651--669.

\bibitem[Mal 72]{Mal72} B.~Malgrange, {\em Equations de Lie I \& II,} J.
Differential Geometry \textbf{6} (1972) 503-522 \& {\bf 7} (1972)
117--141.

\bibitem[Mil 83]{Mil83} J.~Milnor, {\em Remarks on infinite dimensional Lie
groups}. Proceedings of Summer School on quantum Gravity, les
Houches, session XL, North-Holland, 1983.

\bibitem[Olv 96]{Olv96} P.~Olver, {\em Non-associative local Lie
groups}, Journal of Lie Theory \textbf{6} (1996) no. 1, 23--51.

\bibitem[OlvPo 03-1]{OlvPo03-1} P.~Olver, J.~Pohjanpelto, {\em Moving frames for pseudo--groups. I. The Maurer--Cartan forms}, preprint,
University of Minnesota, 2003.

\bibitem[OlvPo 03-2]{OlvPo03-2} P.~Olver, J.~Pohjanpelto, {\em Moving frames for pseudo--groups. II. Differential invariants for submanifolds},
preprint, University of Minnesota, 2003.

\bibitem[RK 97]{RK97}T.~Robart and N.~Kamran, {\em Sur la th\'eorie locale des pseudogroupes de transformations continus infinis,}
Math. Ann. \textbf{ 308} (1997), no. 4, 593--613.

\bibitem[SiSe 65]{SiSe65}I.M.~Singer, S.~Sternberg,{\em The infinite groups of Lie and Cartan. I. The transitive groups,}
J. Analyse Math. \textbf{15} (1965) 1--114.
\end{thebibliography}
\end{document}